\documentclass[reqno]{amsart}
\usepackage{amssymb}
\usepackage{hyperref}
\makeatletter

\DeclareMathOperator{\kod}{Kod}

\def\Bbb{\mathbb}
\def\bcp{\mathbb C\mathbb P}

\newtheorem{main}{Theorem}

\newtheorem{thm}{Theorem}[section]
\newtheorem{prop}[thm]{Proposition}
\newtheorem{cor}[thm]{Corollary}
\newtheorem{lem}[thm]{Lemma}

\newtheorem{defn}[thm]{Definition}

\numberwithin{equation}{section}
\newtheorem{question}[thm]{Question}
\newenvironment{xpl}{\medskip 
\noindent {\bf  Example}}{\mbox{}\medskip}

\begin{document}

	\title{The Kodaira dimension of 
	diffeomorphic K\"ahler 3-folds}
	\author{Rare{\c s} R{\u a}sdeaconu} 
	\address{Department of Mathematics, SUNY at Stony Brook, Stony Brook,
		11794, NY, USA}
	\email{rares@math.sunysb.edu}
	\subjclass[2000]{Primary 32J17; Secondary 19J10}
	\keywords{Whitehead torsion,s-cobordism,Kodaira dimension,
		K{\" a}hler manifolds} 
	\date{\today}

	\begin{abstract}
We provide infinitely many examples of pairs of diffeomorphic, 
non simply connected K\" ahler manifolds of complex dimension 
three with different Kodaira dimensions. Also, in any possible 
Kodaira dimension we find infinitely many pairs of non deformation 
equivalent, diffeomorphic K\" ahler threefolds. 
	\end{abstract}
%\makeatother
\maketitle

	\section{Introduction}
	\label{introduction}

Let $M$ be a compact complex manifold of complex dimension $n.$ On any such 
manifold the {\em canonical line bundle} $K_M=\wedge^{n,0}$ encodes 
important information about the complex structure. One can 
define a series of birational invariants of $M,$ 
$P_n(M) := h^0(M,K^{\otimes n}_M),~ n\geq 0,$ called the {\em plurigenera}.  The 
number of independent holomorphic $n$-forms on $M,~p_g(M)=P_1(M)$ 
is called the geometric genus. The {\em Kodaira dimension} $\kod (M),$ is a 
birational invariant given by:  
$$
\kod (M) = \limsup \frac{\log h^0(M,K^{\otimes n}_M)}{\log n}.
$$
This can be shown to coincide with the maximal complex dimension 
of the image of $M$ under the pluri-canonical maps, so that 
$\kod (M)\in \{ -\infty , 0, 1, \ldots, n\}$. 
A compact complex $n$-manifold is said to be of {\em general type} if 
$\kod (M)=n$. For  Riemann surfaces, the classification with respect to the 
Kodaira dimension, $\kod (M) = -\infty , 0$ or $1$ 
is equivalent to the one given by the {\em genus}, $g(M)=0,~1,$ 
and $\geq 2,$ respectively.

An important question in differential geometry is to understand 
how the complex structures on a given complex manifold are related to the 
diffeomorphism type of the underlying smooth manifold or further, to 
the topological type of the underlying topological manifold. 
Shedding some light on this question is 
S. Donaldson's result on the ``failure of the $h$-cobordism conjecture in 
dimension four''. In this regard, he found a pair of {\em 
non-diffeomorphic}, $h$-cobordant, simply connected 4-manifolds. 
One of them was $\bcp_2\# 9\overline {\bcp}_2,$ the blow-up of $\bcp_2$ 
at nine appropriate points, and the other one was a certain properly elliptic 
surface. For us, an important feature of these two complex surfaces is the fact 
that they have  {\em different Kodaira dimensions}. Later, R. Friedman and Z. Qin 
\cite{frqin} went further and proved that actually, 
for complex surfaces of K\"ahler type, \emph {the Kodaira dimension is 
invariant under diffeomorphisms}. 
However, in higher dimensions, C. LeBrun and F. Catanese gave 
examples \cite{claudecatanese} of pairs of diffeomorphic projective manifolds of 
complex dimensions $2n$ with $n\geq 2,$ and Kodaira dimensions $-\infty$ 
and $2n.$ 
 
In this article we address the question of the invariance of the 
Kodaira dimension under diffeomorphisms in complex dimension $3.$ 
We obtain the expected negative result:
\begin{main}
\label{A}
\label{difkod}
For any allowed pair of distinct Kodaira 
dimensions $(d,d'),$ with the exception of $(-\infty, 0)$ and $(0,3),$ there 
exist infinitely many pairs of diffeomorphic K\"ahler threefolds $(M,M'),$ having 
the same Chern numbers, but with $\kod (M)=d$ and $\kod (M')=d',$ respectively.
\end{main}

\begin{cor}
For K\"ahler threefolds, the Kodaira dimension is not a smooth invariant. 
\end{cor}

Our examples also provide negative answers to questions regarding the 
deformation types of K\"ahler threefolds. 
Recall that two manifolds $X_1$ and $X_2$ are called 
{\em directly deformation equivalent} if there exists a complex manifold 
${\mathcal X},$ and a proper holomorphic submersion 
$\varpi:{\mathcal X}\rightarrow {\Delta}$ with 
$\Delta=\{|z|=1\}\subset \Bbb C,$ such that $X_1$ and $X_2$ occur as 
fibers of $\varpi.$ 
The {\em deformation equivalence} relation is the equivalence relation 
generated by direct deformation equivalence.
 
It is known that two deformation equivalent 
manifolds are orientedly diffeomorphic. For complex surfaces of K\" ahler 
type there were strong indications that the converse should also be true. R. 
Friedman and J. Morgan proved \cite {fm} that, not only the Kodaira dimension is 
a smooth invariant but the plurigenera, too. However, Manetti \cite{manetti} 
exhibited examples of diffeomorphic complex surfaces of general type 
which were {\em not} deformation equivalent. An easy consequence of our $\bf 
{Theorem ~\ref{difkod}}$ and of the deformation invariance of plurigenera for 
3-folds \cite{km} is that in complex dimension 3 the situation is similar:

\begin{cor}
For K\"ahler threefolds the deformation type does not 
coincide with the diffeomorphism type.
\end{cor}
Actually, with a bit more work we can get:
\begin{main}
\label{B}
\label{defdif}
In any possible Kodaira dimension, there exist infinitely many examples of pairs 
of diffeomorphic, non-deformation equivalent K\"ahler threefolds with the same 
Chern numbers.  
\end{main}

The examples we use are Cartesian products of simply connected, $h-$cobordant 
complex surfaces with Riemann surfaces of positive genus. The real six-manifolds 
obtained will therefore be $h-$cobordant. To prove that these six-manifolds are 
in fact diffeomorphic, we use the $s-$Cobordism Theorem, by showing that the 
obstruction to the triviality of the corresponding $h-$cobordism, the {\em 
Whitehead torsion}, vanishes. Similar examples were previously used by Y. Ruan 
\cite{ruan} to find pairs of diffeomorphic symplectic 6-manifolds which are not 
symplectic deformation equivalent. However, to show that his examples are 
diffeomorphic, Ruan uses the classification (up to diffeomorphisms) of {\em 
simply-connected}, 	real 6-manifolds \cite{okonek}. This restricts Ruan's 
construction to the case of Cartesian products by 2-spheres, a result which 
would also follow from Smale's $h$-cobordism theorem.

The examples of pairs complex structures we find are all of K\" ahler type 
with the same Chern numbers. This should be contrasted with 
C. LeBrun's examples  \cite{ctop} of complex structures, 
mostly non-K\" ahler, with 
{\em different~Chern~numbers} on a given differentiable \emph {real} manifold. 

In our opinion, the novelty of this article is the 
use of the apparently forgotten {\em s-Cobordism~Theorem.}  
This theorem is especially useful when combined 
with a theorem on the vanishing of the {\em Whitehead~group}. 
For this, there exist nowadays strong results, due to 
F.T. Farrell and L. Jones \cite{jones}. 

In the next section, we will review the main tools we use to find 
our examples: $h$-cobordisms, the Whitehead group and its vanishing. 
In section 3 we recall few well-known generalities about complex surfaces. 
Sections 4 and 5 contain a number of examples and the proofs of $\bf 
{Theorems~\ref{difkod}}$ and $\bf \ref{defdif}$. In the last section we conclude 
with few remarks and we raise some natural questions.

\section{The s-Cobordism Theorem}
\label{white}
\begin{defn}
\label{defcob}
Let $M$ and $M'$ be two $n$-dimensional closed, smooth, oriented 
manifolds. A $cobordism$ between $M$ and $M'$ is a triplet 
$(W;M,M'),$ where $W$ is an $(n+1)$-dimensional compact, oriented 
manifold with boundary, 
$\partial W=\overline {\partial W}_- \bigsqcup \partial W_+$ with 
${\partial W_-}=M$ and ${\partial W_+}=M'$ 
(by $\overline {\partial W}_-$ we denoted the orientation-reversed 
version of $\partial W_-$). 

We say that the cobordism $(W;M,M')$ is an $h$-cobordism if the 
inclusions $i_-:M\rightarrow W$ and $i_+:M'\rightarrow W$ are 
homotopy equivalences between $M, M'$ and $W.$
\end{defn}
The following well-known results \cite{wall1}, \cite{wall} 
allow us to easily check when two simply connected 4-manifolds are $h$-cobordant:
\begin{thm}
\label{hcobord}
Two simply connected smooth manifolds of dimension 4 are $h$-cobordant 
if and only if their intersection forms are isomorphic.
\end{thm}
\begin{thm}
Any indefinite, unimodular, bilinear form is uniquely determined by 
its rank, signature and parity.
\end{thm}

In higher dimensions any $h$-cobordism $(W;M,M')$ is controlled by a complicated 
torsion invariant $\tau (W;M),$  the {\em Whitehead~torsion,} an element of the 
so called {\em Whitehead~group} which will be defined below.  

Let $\Pi$ be any group, and $R={\Bbb Z}(\Pi)$ the integral unitary 
ring generated by $\Pi.$ We denote by $GL_n (R)$ the group of all 
nonsingular $n\times n$ matrices over $R$. For all $n$ 
we have a natural inclusion $GL_n(R)\subset GL_{n+1}(R)$ identifying 
each $A\in GL_n(R)$ with the matrix: 
$$
\left(\begin{array}{cc} A&0 \\ 0&1
\end{array}
\right)\in GL_{n+1}(R).
$$
Let $\displaystyle GL(R)=\bigcup_{n=1}^{\infty}GL_n (R).$ 
We define the following group: 
$$
K_1(R)=GL(R)/[GL(R),GL(R)].
$$ 
The {\em Whitehead~group} we are interested in is:
$$
Wh(\Pi)=K_1 (R)/<\pm g ~|~ g\in \Pi>.
$$
\begin{thm}
\label{cobordism}
Let $M$ be a smooth, closed manifold. 
For any $h$-cobordism $W$ of $M$ with 
${\partial}_- W=M,$ and with $dim~W\geq 6$  
there exists an element $\tau (W)\in Wh({\pi}_1(M)),$
called {\em the~Whitehead~torsion,}
characterized by the following properties:
	\begin{description}
		\item[s-Cobordism Theorem ]
		\label{cob}
				$\tau (W)\in Wh({\pi}_1(M))=0$ 
				if and only if the \\ 
				\indent $h$-cobordism is trivial, i.e.
				$W$ is diffeomorphic to ${\partial}_- W\times [0,1];$ 
		\item[Existence ]Given $\alpha \in Wh({\pi}_1(M)),$ there exists an 
				$h$-cobordism $W$ with \\ 
				\indent $\tau (W)=\alpha;$
		\item[Uniqueness ]$\tau (W)=\tau (W')$ if and only if there 
				exists a diffeomorphism \\
				\indent $h:W\rightarrow W'$ such that $h_{|M}=id_{M}.$
	\end{description} 
\end{thm}

For the definition of the $Whitehead~torsion$ and the above theorem we refer the 
reader to Milnor's article \cite{milnor}. However, the above theorem 
suffices. When $M$ is simply connected, the s-cobordism theorem is 
nothing but the usual $h$-cobordism theorem \cite{jack}, due to Smale.

This theorem will be a stepping stone in finding pairs 
of diffeomorphic manifolds in dimensions greater than 5, 
provided knowledge about the vanishing of the $Whitehead$ $groups.$ 
The most powerful vanishing theorem that we are aware of 
is the following:
\begin{thm}[Farrell, Jones]
\label{vanishing}
Let $M$ be a compact Riemannian manifold of non-positive sectional 
curvature. Then $Wh({\pi}_1 (M))=0.$
\end{thm}
The uniformization theorem of compact Riemann surfaces yields then the following 
result which, as it was kindly pointed to us by L. Jones, was also known to F. 
Waldhausen \cite{wald}, long before \cite{jones}.
\begin{cor}
\label{rsurf}
Let $\Sigma$ be a compact Riemann surface. Then $Wh({\pi}_1 (\Sigma))=0.$
\end{cor}
An useful corollary, which will be frequently used is the following:  
\begin{cor}
\label{main}
Let $M$ and $M'$ be two simply connected, $h$-cobordant 4-manifolds, 
and $\Sigma$ be a Riemann surface of positive genus. Then $M\times \Sigma$ 
and $M'\times \Sigma$ are diffeomorphic. 
\end{cor}
\proof
Let $W$ be an $h$-cobordism between $M$ and $M'$ such that 
${\partial}_- W=M$ and ${\partial}_+ W=M'$ and let 
${\tilde W}=W\times \Sigma.$ Then 
${\partial}_- \tilde W =M\times \Sigma,~{\partial}_+\tilde W=M'\times 
\Sigma,$ and $\tilde W$ is an $h$-cobordism between 
$M\times \Sigma$ and $M'\times \Sigma.$ Now, since M is simply connected $\pi_1 
(M\times \Sigma)=\pi_1(\Sigma)$ and so 
$Wh(\pi_1 (M\times \Sigma))=Wh(\pi_1(\Sigma)).$ By the 
uniformization theorem any Riemann surface of positive genus admits a metric of 
non-positive curvature. Thus, by Theorem \ref{vanishing},  
$Wh(\pi_1(\Sigma))=0,$ which, by Theorem \ref{cob}.1, implies that 
$M\times \Sigma$ and $M'\times \Sigma$ are diffeomorphic.   
\qed

\section{Generalities}
\label{gen}

To prove $\bf {Theorems ~\ref{difkod}}$ and 
$\bf \ref{defdif}$ we will use our Corollary 
\ref{main}, by taking for $M$ and $M'$ appropriate $h$-cobordant, 
simply connected, complex projective surfaces, and for $\Sigma$, 
Riemann surfaces of genus $g(\Sigma)\geq 1.$ 
To find examples of $h$-cobordant complex surfaces, we use:

\begin{prop}
\label{cobcrit}   
Let $M$ and $M'$ be two simply connected complex surfaces 
with the same geometric genus $p_g,~c^2_1(M)-c^2_1(M')=m\geq 0$ 
and let $k>0$ be any integer. Let $X$ be the blowing-up of $M$ 
at $k+m$ distinct points and $X'$ be the blowing-up of $M'$ at $k$ 
distinct points. Then $X$ and $X'$ are $h$-cobordant, 
$\kod (X)=\kod (M)$ and $\kod (X')=\kod (M').$
\end {prop}
\proof 
By Noether's formula we see that $b_2(M')=b_2(M)+m.$ Since, by 
blowing-up we increase each time the second Betti number by one, 
it follows  that $b_2(X')=b_2(X).$ 
Using the birational invariance of the plurigenera,  
we have that $b_+(X')=2p_g+1=b_+(X).$ As $X$ and $X'$ are both 
non-spin, and their intersection forms have the same rank and signature, their 
intersection forms are isomorphic.Thus, by Theorem \ref{hcobord}, $X$ and 
$X'$ are $h$-cobordant. The statement about the Kodaira dimension follows from 
its birational invariance, too.     
\qed 
%
%Using the fact that by blowing up, the Kodaira dimension doesn't change, 
%an easy corollary of the proof of Proposition \ref{cobcrit} is the following:
\begin{cor}
\label{infinity}
Let $S$ and $S'$ be two simply connected, $h$-cobordant complex surfaces.
If $S_k$ and $S_k'$ are the blowing-ups of the two surfaces, each at $k\geq 0$ 
distinct points, then $S_k$ and $S_k'$ are $h$-cobordant, too. Moreover, 
$\kod (S_k)=\kod (S),$ and  $\kod (S_k')=\kod (S').$
\end{cor}
The following proposition will take care of the computation of the 
Kodaira dimension of our examples. Its proof is standard, 
and we will omit it.
\begin{prop}
\label{kodim}
Let $V$ and $W$ be two complex manifolds. Then 
$P_m(V\times W)=P_m(V)\cdot P_m(W).$ In particular,  
$\kod (V\times W)=\kod (V)+\kod (W).$
\end{prop}
For the computation of the Chern numbers of the examples involved, 
we need:
\begin{prop}
\label{chern}
Let $M$ be a smooth complex surface with $c^2_1(M)=a,~c_2(M)=b,$ 
and let $\Sigma$ be a smooth complex curve of genus $g,$ and $X=M\times \Sigma$ 
their Cartesian product. The Chern numbers $(\bf {c^3_1,c_1c_2,c_3})$ of $X$ are 
$((6-6g)a,(2-2g)(a+b),(2-2g)b).$
\end{prop}
\proof
Let $p:X\rightarrow M,$ and $q:X\rightarrow {\Sigma}$ 
be the projections onto the two factors. Then the total Chern class 
is: $c(X)=p^*c(M)\cdot q^*c(\Sigma),$ which allows us to identify the Chern 
classes. Integrating over $X,$ the result follows immediately.
\qed 

\section{Diffeomorphism types - Proof of Theorem \ref{A}}
\label{examples}

In this section we prove  $\bf Theorem ~\ref{difkod}$. To do this, for each of 
the pairs of Kodaira dimensions stated, we provide infinitely many examples, by 
taking Cartesian products of appropriate $h-$cobordant 
K\"ahler surfaces with Riemann surfaces of positive genus.  

\begin{xpl}{~\bf 1:}{~\em Pairs of Kodaira dimensions $(-\infty,1)$ and 
$(-\infty,2)$}
\label{dolgachev}

\smallskip
Let $M$ be the blowing-up of $\bcp_2$ at $9$ distinct points given by the 
intersection 
of two generic cubics. $M$ is a non-spin, simply connected complex surface 
with $\kod(M)=-\infty$ which is also an elliptic fibration,  
$\pi :M\rightarrow \bcp_1.$ By taking the cubics general enough, we may assume 
that $M$ has no multiple fibers, and the only singular fibers are irreducible 
curves with one ordinary double point. Let $M'$ be obtained from $M$ by 
performing logarithmic transformations on two of its smooth fibers, with  
multiplicities $p$ and $q,$ where $p$ and $q$ are two relatively prime positive 
integers. $M'$ is also an elliptic surface, ${\pi}':M'\rightarrow \bcp_1,$ whose 
fibers can be identified to those in $M$ except for the pair of multiple fibers 
$F_1,$ and $F_2.$ Let $F$ be homology class of the generic fiber in $M'.$
In homology we have $[F]=p[F_1]=q[F_2].$
By canonical bundle formula, we see that:
$K_M=-F,$ and 
\begin{equation}
\label{cbf}
K_{M'}=-F+(p-1)F_1+(q-1)F_2=\frac{pq-p-q}{pq}F.
\end{equation}
Then $p_g(M)=p_g(M')=0,~c^2_1(M)=c^2_1(M')=0,~\text{and}~\kod (M')=1.$
Moreover, from \cite[Theorem 2.3, page 158]{frmorgan} $M'$ is 
simply connected and non-spin. 
 
For any $k\geq 0,$ let $M_k$ and $M_k'$ be the blowing-ups at $k$ 
distinct points of $M$ and $M',$ 
respectively, and let $ \Sigma$ be a Riemann surface. If 
$g(\Sigma)=1,$  according to 
Corollary \ref{main} and Proposition \ref{kodim},
$(M_k\times \Sigma_1,M_k'\times \Sigma_1),~k\geq 0$ will provide 
infinitely many pairs of diffeomorphic K{\"a}hler threefolds, of Kodaira 
dimensions are $-\infty$ and $1,$ respectively. If $g(\Sigma)\geq2,$ 
we get infinitely many pairs of 
diffeomorphic K\"ahler threefolds with Kodaira dimensions are $-\infty,$ 
and $2,$ respectively. The statement about the Chern numbers follows from 
Proposition \ref{chern}.
\end{xpl}

\smallskip
\begin{xpl}{~\bf 2:}{~\em Pairs of Kodaira dimensions $(0,1)$ and $(0,2)$}
\label{homotpyk3}

\smallskip
In $\bcp_1 \times \bcp_2,$ let $M$ be the the generic section of
line bundle 
$p_1^*{\mathcal O}_{\bcp_1}(2)\otimes p_2^*{\mathcal O}_{\bcp_2}(3),$ 
where $p_i,$ $i=1,2$ are the projections onto the two factors.
Then $M$ is a $K3$ surface, i.e. a smooth, simply connected complex surface, 
with trivial canonical bundle. Moreover,   
using the projection onto the first factor, it fibers over $\bcp_1$ 
with elliptic fibers. 

Kodaira \cite{kodaira} produced infinitely many examples of properly elliptic 
surfaces of K\" ahler type, homotopically equivalent to a $K3$ surface, 
by performing two logarithmic transformations on two smooth fibers with 
relatively prime multiplicities on such elliptic $K3.$ Let $M'$ to be any 
such surface, and let $M_k$ and $M_k'$ be the blowing-ups at 
$k$ distinct points of $M$ and $M',$ respectively. As before, let $\Sigma$ be 
a Riemann surface. If $g(\Sigma)=1,$ the Cartesian products 
$M_k\times \Sigma$ and  $M_k'\times \Sigma$ will provide infinitely many 
pairs of diffeomorphic  K\"ahler $3$-folds of Kodaira 
dimensions $0$ and $1,$ respectively. If $g(\Sigma)\geq 2,$  
we obtain pairs in Kodaira 
dimensions $1$ and $2,$ respectively. Again, the statement about the 
Chern numbers follows from Proposition \ref{chern}.
\end{xpl}

\smallskip
\begin{xpl}{~\bf 3:}{~\em Pairs of Kodaira dimensions 
$(-\infty,2)$ and $(-\infty,3)$}
\label{bar}

\smallskip
Arguing as before, we present a pair of 
simply connected, $h-$cobordant projective surfaces, one on Kodaira 
dimension 2, and the other one of Kodaira dimension $-\infty.$

Let $M$ be the {\em Barlow surface} \cite{barlow}. This is a non-spin, 
simply connected projective surface of general type, with $p_g=0$ and 
$c^2_1(M)=1.$ 
It is therefore $h$-cobordant to $M'$, the projective plane ${\bcp}_2$ 
blown-up at 8 points. 
By taking the Cartesian product of their blowing-ups by 
a Riemann surface of genus 1, we obtain diffeomorphic, 
projective threefolds of Kodaira dimensions $3,$ and $-\infty,$ respectively, 
while for a Riemann surface of bigger genus, we obtain  diffeomorphic, 
projective threefolds of Kodaira dimensions $2,$ and $-\infty,$ respectively. 
The invariance of their Chern numbers follows as usual.
\end{xpl}

\begin{xpl}{~\bf 4:}{~\em  Pairs of Kodaira dimensions $(0,2)$ and $(1,3)$}
\label{catanese}

\smallskip
Following \cite{catanese}, we will describe an example of simply connected, 
minimal surface of general type with $c^2_1=p_g=1.$

In $\bcp_2$ we consider two generic smooth cubics $F_1$ and $F_2$, 
which meet transversally at 9 distinct points, $x_1,\cdots,x_9,$ 
and let $\sigma :\tilde X\rightarrow \bcp_2$ be the blowing-up of $\bcp_2$ 
at $x_1,\cdots,x_9,$ with exceptional divisors ${\tilde E}_i,~ i=1,...,9.$ 
Let ${\tilde F}_1$ and ${\tilde F}_2$ be the strict transforms of 
$F_1$ and $F_2,$ 
respectively. Then ${\tilde F}_1$ and ${\tilde F}_2$ are two disjoint, 
smooth divisors, and we can easily see that 
${\mathcal O}_{\tilde X}({\tilde F}_1+{\tilde F}_2)=
{\tilde {\mathcal L}}^{\otimes 2},$ 
where 
$$
{\tilde {\mathcal L}}={\sigma}^*{\mathcal O}_{\bcp_2}(3)\otimes
{\mathcal O}_{\tilde X}({\tilde E}_1+\dots+{\tilde E}_9).
$$

Let $\pi:{\bar X}\rightarrow {\tilde X}$ to be the 
double covering of ${\tilde X}$ branched along the smooth divisor 
${\tilde F}_1+{\tilde F}_2.$ We denote by $p:{\bar X}\rightarrow {\bcp_2}$ the 
composition $\sigma \circ \pi,$ and by ${\bar F_1},~{\bar F_2}$ the reduced 
divisors ${\pi}^{-1}({\tilde F_1}),$ and ${\pi}^{-1}({\tilde F_2}),$ respectively.
Since each $\tilde E_i$ intersects the branch locus at 2 distinct 
points, we can see 
that for each $i=1,\dots,9,~\bar E_i={\pi}^{-1}({\tilde E_i})$ 
is a smooth (-2)-curve  
such that ${\pi}_{|{\bar E_i}}:{\bar E_i}\rightarrow {\tilde E_i}$ is the double 
covering of $\tilde E_i$ branched at the two intersection 
points of $\tilde E_1$ with 
${\tilde F}_1+{\tilde F}_2.$ As the $\tilde E_i~'$s are mutually disjoint, 
the $\bar E_i~'$s 
will also be mutually disjoint. Similarly, if $\ell$ is a line in 
$\bcp_2$ not passing through any of the intersection points of $F_1$ 
with $F_2,$ then  $L=p^*({\ell})=p^*{\mathcal O}_{\bcp_2}(1)$ 
is a smooth curve of genus 2, not intersecting any of the $\bar E_i~'$s. 
Since 
$$
p^*{\mathcal O}_{\bcp_2}(3)={\mathcal O}_{\bar X}(2{\bar F_1}+ 
{\bar E_1}+\dots+{\bar E_9}),
$$
we can write as before ${\mathcal O}_{\bar X}(L+{\bar E_1}+\dots+{\bar E_9})
={\bar {\mathcal L}}^{\otimes 2},$
where 
$$
{\bar {\mathcal L}}=p^*{\mathcal O}_{\bcp_2}(2)\otimes 
{\mathcal O}_{\bar X} (-{\bar F}_1).
$$

Let now $\phi:{\bar S}\rightarrow {\bar X}$ be the double covering of 
$\bar X$ ramified along the smooth divisor $L+{\bar E_1}+\dots+{\bar 
E_9}.$ The surface $\bar S$ is non-minimal with exactly 9 disjoint exceptional 
curves of the first kind. Namely, the reduced divisors 
${\phi}^{-1}({\bar E_i}),~i=1,\dots 9.$ The surface $S$ we were 
looking for is obtained from $\bar S$ by blowing down these 9 exceptional 
curves. 

F. Catanese proves \cite{catanese} that S is a simply connected, 
minimal surface of general type with $c^2_1(S)=p_g(S)=1.$

Let $S_k'$ be the blowing-up of a $K3$ surface at $k$ distinct points. 
Let also $S_k$ denote the blowing-up of $S$ at $k+1$ distinct points, and 
let $\Sigma$ be a Riemann surface. If 
$g(\Sigma)=1,~(S_k\times \Sigma ,S_k'\times \Sigma)$ will provide 
infinitely many pairs of diffeomorphic 
K\" ahler threefolds of Kodaira dimensions $2$ and $0,$ respectively, while if
$g(\Sigma)\geq 2$ we get 
infinitely many pairs of diffeomorphic 
K\" ahler threefolds of Kodaira dimensions $3$ and $1,$ respectively. 
The statement about the Chern classes follows as before.  
\end{xpl}

\begin{xpl}{~\bf 5:}{~\em Pairs of Kodaira dimensions $(1,2)$ and $(2,3)$}
\label{elliptic}

\smallskip
%\noindent
In $\bcp_1 \times \bcp_2,$ let $M_n$ be the the generic section of
line bundle 
$p_1^*{\mathcal O}_{\bcp_1}(n)\otimes p_2^*{\mathcal O}_{\bcp_2}(3)$ 
for $n\geq 3,$ where $p_i,$ $i=1,2$ be the projections onto the two factors.
Then $M_n$ is a smooth, simply connected projective surface, and 
using the projection onto the first factor we see that $M_n$ is a properly 
elliptic surface. By the adjunction formula, 
the canonical line bundle is:
$$K_{M_n}=p_1^*{\mathcal O}_{\bcp_1}(n-2).$$ 

From this and the projection formula we can find the purigenera: 
	\begin{equation}
		\begin{aligned}
P_m(M_n)&=h^0(M_n,K_{M_n}^{\otimes m})=h^0(M_n,p_1^*{\mathcal 
O}_{\bcp_1}(m(n-2)))\\
&=h^0(\bcp_1,{p_1}_*p_1^*{\mathcal O}_{\bcp_1}(m(n-2)))\\
&=h^0(\bcp_1,{\mathcal O}_{\bcp_1}(m(n-2)))\\
&=m(n-2)+1. \notag
		\end{aligned}
	\end{equation}
So, $\kod (M_n)=1,$ and $p_g(M_n)=n-1.$ We can also see that $c^2_1(M_n)=0.$ 

Let $M'$ be any smooth sextic in ${\bcp}_3.~M'$ is a simply connected 
surface of general type with $p_g(M')=10,$ and $c^2_1 (M')=24.$ Let $M_k'$ 
be the blowing-up of $M$ at $24+k$ 
distinct points, $M_k$ be the blowing-up of $M_{11}$ at $k+1$ points, 
and let $\Sigma$ be a Riemann surface. If 
$g(\Sigma)=1,~(M_k\times \Sigma ,M_k'\times \Sigma)$ will provide 
infinitely many pairs of diffeomorphic 
K\" ahler threefolds of Kodaira dimensions $1$ and $2,$ respectively, while if
$g(\Sigma)\geq 2$ we get 
infinitely many pairs of diffeomorphic 
K\" ahler threefolds of Kodaira dimensions $2$ and $3,$ respectively. 
The statement about the Chern classes again follows.    
\end{xpl}

\section{Deformation Type - Proof of Theorem \ref{B}}
\label{deftype}

Similar idea can be used to prove $\bf {Theorem ~\ref{defdif}}.$ The proof 
follows from the examples below. 

\noindent
\begin{xpl}{~\bf 1:}{~\em Kodaira dimension $-\infty$}

\smallskip
Here we use again the Barlow surface $M,$ and $M',$ 
the blowing-up of $\bcp_2$ at 8 points as two $h$-cobordant complex surfaces. Let 
$S_k$ and $S_k'$ denote the blowing-ups of $M$ and $M',$ respectively at $k$ 
distinct points. Then, by the classical $h$-cobordism 
theorem, $X_k=S_k\times \bcp_1$ and $X_k'=S_k'\times \bcp_1$ are two 
diffeomorphic $3$-folds with the same Kodaira dimension $-\infty.$ The fact that 
$X_k$ and $X_k'$ are not deformation equivalent follows as in 
\cite{ruan} from Kodaira's stability 
theorem \cite{kstab}. We also see immediately that they have the same Chern 
numbers.   
\end{xpl}

\smallskip
\noindent
\begin{xpl}{~\bf 2:}{~\em Kodaira dimension $2$ and $3$}

\smallskip
We start with a {\em Horikawa surface}, namely a simply connected surface of 
general type $M$ with $c^2_1(M)=16$ and $p_g(M)=10.$ An example of such 
surface can be obtained as a ramified double cover of $Y=\bcp_1\times 
\bcp_1$ branched at a generic curve of bi-degree $(6,12).$ If we denote by 
$p:M\rightarrow Y,$ its degree $2$ morphism onto $Y,$ then the canonical 
bundle of $M$ is $K_M={\mathcal O}_Y(1,4),$ see \cite[page 182]{barth}. 
Here by ${\mathcal O}_Y(a,b)$ we denote the line bundle 
$p^*_1{\mathcal O}_{\bcp_1}(a)\otimes p^*_2{\mathcal O}_{\bcp_1}(b),$ 
where $p_i,~1=1,2$ are the projections of $Y$ onto the two factors. 
Notice that the formula for the canonical bundle shows that $M$ is not spin.
\begin{lem}
\label{hor}
The plurigenera of $M$ are given by: 
$$
P_n(M)=\left\{ \begin{array}{cc}
10&n=1 \\
8n^2-8n+11&n\geq2 \end{array} \right.
$$
\end{lem}
\proof
Cf. \cite{barth} we have 
$p_*{\mathcal O}_M={\mathcal O}_Y\oplus {\mathcal O}_Y(-3,-6).$ 
We have:
	\begin{align*}
P_n(M)&=h^0(M,p^*{\mathcal O}_Y(n,4n))= 
h^0(Y,p_*p^*{\mathcal O}_Y(n,4n))\\ 
&=h^0(Y,{\mathcal O}_Y(n,4n)\otimes p_*{\mathcal O}_M)\\ 
&=h^0(Y,{\mathcal O}_Y(n,4n))+h^0(Y,{\mathcal O}_Y(n-3,4n-6)).
	\end{align*}
Now, if $n<3$ we get 
$P_n(M)=(n+1)(4n+1).$ In particular, $p_g(M)=10$ and $P_2(M)=27.$  
If $n\geq 3,~P_n(M)=(n+1)(4n+1)+(n-2)(4n-5)=8n^2-8n+11.$ 
\qed

\smallskip
Let $M' \subset \bcp_3$ be a smooth sextic. The adjunction formula will provide 
again the the canonical bundle $K_{M'}={\mathcal O}_{M'}(2)$ and so 
$c^2_1(M')=24.$ 
\begin{lem}
\label{sextic}
The plurigenera of $M'$ are given by: 
$$
P_n(M')=\left\{ \begin{array}{cc} \binom{2n+3}{3}&{n=1,2} \\
{12n^2-12n+11}&{n\geq 3} \end{array} \right.
$$
\end{lem} 
\proof
From the exact sequence 
$0\rightarrow {\mathcal O}_{\bcp_3}(2n-6)\rightarrow 
{\mathcal O}_{\bcp_3}(2n)\rightarrow K^{\otimes n}_{M'}\rightarrow 0,$
we get: 
		\begin{align*} 
0&\rightarrow H^0(\bcp_3,{\mathcal O}_{\bcp_3}(2n-6))\rightarrow  
H^0(\bcp_3,{\mathcal O}_{\bcp_3}(2n)) \\
&\rightarrow H^0({M'},K^{\otimes n}_{M'})
\rightarrow H^1(\bcp_3,{\mathcal O}_{\bcp_3}(2n))=0. 
		\end{align*}
So, for $n\geq3,$ 
$$
P_n(M')=\binom{2n+3}{3}-\binom{2n-3}{3}=12n^2-12n+11,
$$
while for $n<3,~P_n(M')=\binom{2n+3}{3}.$ In particular, 
$p_g(M')=10$ and $P_2(M')=35.$ 
\qed

Let $M_k$ be the blowing-up of $M$ at $k$ distinct points, $M_k'$ be the 
blowing-up of $M'$ at $8+k$ distinct points, and let $\Sigma$ be a Riemann 
surface. If $g(\Sigma)=1,~(M_k\times \Sigma , M_k'\times \Sigma),k\geq 0$ 
will provide the required examples of Kodaira dimension 2, and if 
$g(\Sigma)\geq 2,$ will provide the required examples of Kodaira dimension 3. 

To prove that they are not deformation 
equivalent we will use the deformation invariance of plurigenera theorem 
\cite[page 535]{km}. Because of the their multiplicative property cf. Proposition 
\ref{kodim}, it will suffice to look at the plurigenera of $M$ and $M'.$ But, 
$P_2(M)=27$ and $P_2(M')=P_2(S)=35,$ and so $M\times \Sigma$ and $M'\times 
\Sigma$ are not deformation equivalent. 

The statement about the Chern numbers of this examples follows immediately. 
\end{xpl}

\noindent
\begin{xpl}{~\bf 3:}{~\em Kodaira dimension $1$}

\smallskip
Here we use again the elliptic 
surfaces $\pi:M_{p,q}\rightarrow {\bcp_1}$ obtained from the rational elliptic 
surface by applying logarithmic transformations 
on two smooth fibers, with relatively prime multiplicities $p$ and $q.$ 
%These surfaces are non-spin, simply connected and their canonical 
%bundle is $K_{M_{p,q}}=-F+(p-1)F_1+(q-1)F_2=\frac{pq-p-q}{pq}F,$
%where F is the class of a generic fiber. 
From (\ref{cbf}) we get 
$K^{\otimes pq}_{M_{p,q}}=p^*{\mathcal O}_{\bcp_1}((pq-p-q)).$ Hence 
$P_{pq}(M_{p,q})=pq-p-q+1,~{\text {while if}} ~n\leq pq, 
~P_n(M_{p,q})=0,$ 
the class of $F$ being a primitive element in $H^2(M_{p,q},\Bbb Z),$ cf. 
\cite{kodaira}. It is easy to see now that, for example, if 
$(p,q)\neq (2,3),~P_6(M_{p,q})\neq P_6(M_{2,3}).$ If $\Sigma$ is any 
smooth elliptic curve, the $3-$folds  
$X_{p,q}=M_{p,q}\times \Sigma$ will provide infinitely many diffeomorphic 
K\"ahler threefolds of Kodaira dimension $1.$ Corollary \ref{chern} 
shows again that all these threefolds have the same Chern numbers.
The above computation of plurigenera shows that, in general, 
the $X_{p,q}$'s have different plurigenera. Hence, 
these K\"ahler threefolds are not deformation equivalent. 
\end{xpl}

\noindent
\begin{xpl}{~\bf 4:}{~\em Kodaira dimension $0$}

\smallskip
Here we are supposed to start with a simply connected minimal 
surface of zero Kodaira dimension. But, up to diffeomorphisms there 
exists only one \cite{barth}, the $K3$ surface. So our method fails 
to produce examples in this case.
However, M. Gross constructed \cite{gross} a pair of diffeomorphic 
complex threefolds with trivial canonical bundle, which are {\em not} 
deformation equivalent. For the sake of completeness we will briefly 
recall his examples. 

Let $E_1={\mathcal O}^{\oplus 4}_{\bcp_1}$ and 
$E_2={\mathcal O}_{\bcp_1}(-1)\oplus {\mathcal O}^{\oplus 2}_{\bcp_1}(1)
\oplus {\mathcal O}_{\bcp_1}$ be the two rank 4 vector bundles over 
$\bcp_1,$
and consider $X_1={\Bbb P}(E_1)$ and $X_2={\Bbb P}(E_2)$ their  
projectivizations. Note that $E_2$ deforms to $E_1.$ Let  
$M_i\in |-K_{X_i}|,~i=1,2$ general anticanonical divisors. 
The adjunction formula immediately shows that $K_{M_i}=0,~i=1,2,$ and so 
$M_1$ and $M_2$ have zero Kodaira dimension. 
While for $M_1$ is easy to see that can be chosen to be smooth, simply 
connected and with no torsion in cohomology, Gross shows \cite{gross}, 
\cite{ruangross} that the same holds for $M_2.$ 
Moreover, the two 3-folds have the same topological invariants, 
(the second cohomology group, the Euler characteristic, the cubic form, 
and the first Pontrjaghin class), and so, cf. \cite{okonek}, 
are diffeomorphic. To show that $M_1$ and $M_2,$ 
are not deformation equivalent, note that $M_2$ 
contains a smooth rational curve with normal bundle 
${\mathcal O}(-1)\oplus {\mathcal O}(-1),$ which is stable under the 
deformation of the complex structure while $M_1,$ doesn't. 
Obviously, $M_1$ and $M_2$ have the same Chern numbers. By 
blowing them up simultaneously at 
$k$ distinct points, we obtain infinitely many pairs of diffeomorphic, 
projective threefolds of zero Kodaira dimension with 
the same Chern numbers.
\end{xpl}
\section{Concluding Remarks}
\label{cr}
$\bf 1.$ ~Let $M$ and $M'$ be any of the pairs of complex surfaces 
discussed in the previous two sections. A simple inspection shows that 
they are not spin, and so, their intersection forms will have the form 
$m\langle 1\rangle \oplus n\langle -1\rangle.$ By a result of Wall 
\cite{wall1}, if $m,n\geq 2,$ the intersection form is transitive on the 
primitive characteristic elements of fixed square.  Since, $c_1$ is 
characteristic, if it is primitive too, we can assume that the homotopy 
equivalence $f:M\rightarrow M'$ given by an automorphism 
of such intersection form will carry the first Chern class of $M'$ into 
the first Chern class of $M.$ But this implies that the $h-$cobordism 
constructed between $X=M\times {\Sigma}$ and $X'=M'\times {\Sigma}$ also 
preserves the first Chern classes. 

Following Ruan \cite{ruan}, we can arrange our examples such that $c_1$ 
is a primitive class. In the cases when $b_+>1,$ which is equivalent to 
$p_g>0,$ it follows that there exists a diffeomorphism $F:X\rightarrow X'$ 
such that $F^*c_1(X')=c_1(X),$ where $F^*:H^2(X',\Bbb Z)\rightarrow  
H^2(X,\Bbb Z)$ is the isomorphism induced by $F.$ In these cases our 
theorems provide either examples of pairs of diffeomorphic K\"ahler threefolds, 
with the same Chern classes, but with different Kodaira dimensions, 
or examples of pairs of non deformation equivalent, diffeomorphic K\"ahler 
threefolds, with same Chern classes and of the same Kodaira dimension. 

However, in some cases we are forced to consider surfaces with $b_+=1.$ 
In these cases it is not clear whether one can arrange the $h-$cobordisms 
constructed between $X=M\times {\Sigma}$ and $X'=M'\times {\Sigma}$ also 
preserves the first Chern classes.  

\smallskip
	\indent $\bf 2.$ ~With our method it is impossible to provide examples 
of diffeomorphic $3$-folds of Kodaira dimensions $(0,3)$ and 
$(-\infty, 0).$ In the first case, our method fails for obvious 
reasons. In the second case, the reason is that 
for a projective surface of Kodaira dimension $-\infty,$ the geometric 
genus $p_g$ is $0,$ while for a simply connected projective surface of 
Kodaira dimension $0,$ $p_g\neq 0.$ Thus, any two surfaces of these 
dimensions will have different $b_+,$ which is preserved under blow-ups. 
So, no pair of projective surfaces of these Kodaira dimensions can 
be $h$-cobordant. However, this raises the following question:

\begin{question}
Are there examples of pairs of diffeomorphic, projective $3$-folds 
$(M,M')$ of Kodaira dimensions $(0,3)$ or $(-\infty, 0)?$
\end{question}

Most of the examples exhibited here have the fundamental group of a 
Riemann surface. Natural questions to ask would be the following:

\begin{question}
Are there examples of diffeomorphic, {\em simply connected}, complex, 
projective $3-$folds of different Kodaira dimension?  
\end{question}
\begin{question}
Are there examples of projective, {\em simply connected}, diffeomorphic, but 
not deformation equivalent $3$-manifolds with the same Kodaira dimension?  
\end{question}
As we showed, the answer is {\emph {yes}} when the Kodaira 
dimensions is $-\infty$ or 0, but we are not aware of such examples in the other 
cases.
\section*{Acknowledgements} The author would like to thank his thesis 
advisor Claude LeBrun for suggesting these problems and Lowell Jones for 
some illuminating discussions about the Whitehead groups.

  \end{document}